\documentclass[12pt,a4paper]{amsart}
\usepackage{graphicx,epsfig}
\usepackage{amssymb}
\usepackage{todonotes}
\usepackage[polish,english]{babel}
\usepackage{polski}

\newenvironment{namelist}[1]{%
\begin{list}{}
     {
      
      \settowidth{\labelwidth}{#1}
      \setlength{\leftmargin}{1.1\labelwidth}
               }
      }{%
\end{list}}

\newtheorem{thm}{Theorem}
\newtheorem{lem}{Lemma}
\newtheorem{prop}{Proposition}
\newtheorem{cor}{Corollary}
\newtheorem{Ex}{Example}[section]
\newtheorem{defn}{Definition}
\newtheorem{rem}{Remark}
\begin{document}
\title[Renewal theory for Kendall random walks]{ Renewal theory for extremal Markov sequences of Kendall type}
\author[B.H. Jasiulis-Go{\l}dyn, J.K. Misiewicz, K. Naskr\k{e}t, E. Omey]{Barbara H. Jasiulis-Go{\l}dyn$^1$, Jolanta K. Misiewicz $^2$, Karolina Naskr\k{e}t$^3$ and Edward Omey$^4$}
\thanks{$^1$ Institute of Mathematics, University of Wroc{\l}aw, pl. Grunwaldzki 2/4, 50-384 Wroc{\l}aw, Poland, e-mail: jasiulis@math.uni.wroc.pl \\
$^2$ Faculty of Mathematics and Information Science, Warsaw
University of Technology, ul. Koszykowa 75, 00-662 Warsaw, Poland, e-mail:
j.misiewicz@mini.pw.edu.pl \\
$^3$ Institute of Mathematics, University of Wroc{\l}aw, pl. Grunwaldzki 2/4, 50-384 Wroc{\l}aw, Poland, e-mail: karolina.lukaszewicz@math.uni.wroc.pl\\
$^4$ Faculty of Economics and Business-Campus Brussels, KU Leuven, Warmoesberg 26, 1000 Brussels, Belgium, e-mail:
edward.omey@kuleuven.be 
}
\date{\today}
\begin{abstract}
The paper deals with renewal theory for a class of extremal Markov sequences connected with the Kendall convolution. We consider here some particular cases of the Wold processes associated with generalized convolutions. We prove an analogue of the Fredholm theorem for all regular generalized convolutions algebras. Using regularly varying functions we prove a Blackwell theorem for renewal processes defined by Kendall random walks. 
\end{abstract}
\maketitle

\tableofcontents
\section{Introduction}
In this paper we study renewal theory for Kendall random walks as defined in \cite{KendallWalk}, \cite{factor}. These random walks form a class of extremal Markov sequences (\cite{Alpuim}) with transition probabilities given by a generalized convolution, called the Kendall convolution. Their structure is similar to Pareto processes (\cite{Arnold1}, \cite{Ferreira}), minification  processes (\cite{Lewis}, \cite{Lopez}), the max-autoregressive moving average processes MARMA (\cite{Ferreira}) and extremal processes (\cite{Embrechts}). 

\vspace{3mm}
The origin of generalized convolutions comes from the paper of Kingman (see \cite{King}) about spherically symmetric random walks. Instead of the classical convolution  corresponding to a sum of independent random elements, we consider commutative and associative binary operations described by Urbanik in a series of papers (\cite{Urbanik64}) and further developed by Bingham (\cite{Bin71}-\cite{Bin84}) in the context of regularly varying functions. 
The result of the generalized convolution of two trivial probability measures can be a nondegenerate probability measure. In this paper we focus on the Kendall convolution case. In the Kendall convolution, the result of two probability measures concentrated at 1 is the Pareto distribution with density $\pi_{2\alpha}(dy) =2\alpha y^{-2\alpha-1} \mathbf{1}_{[1,\infty)}(y)dy$. It is the main reason why it produces heavy tailed distributions. Theory of stochastic processes under generalized convolutions was introduced in \cite{BJMR} and following this paper we develop here renewal theory for a class of Markov chains generated by generalized convolutions. In many results we focus on the Kendall convolution case.

\vspace{3mm}

\vspace{3mm}
An ordinary renewal process deals with successive occurences of events
such as failures of machines, arrivals in a queue, lifetimes of systems and
so on. We assume that successive waiting times are given by $T_{1},T_{2},...$
where the $T_{i}$ are independent random variables with a common
distribution function $F$. The time of the $n-$th occurence is given by the
sum $S_{n}=T_{1}+T_{2}+...+T_{n}$ with the convention $S_{0}=0$. If $F(0)=0$, the sequence $(S_{n})$ is called a renewal sequence. If the event $\left\{ S_{n}\leq t\right\} $ is called a "success", then the
total number of successes in the sequence equals 
$$
R(t)=\sum_{n=0}^{\infty }P(S_{n}\leq t)=\sum_{n=0}^{\infty }F^{\ast n}(t).
$$
The function $R(t)$ is called the renewal function and it is well-defined
for all $t$ with $F(t)<1$.

The elementary renewal theorem states that 
\begin{equation}
\lim_{t\rightarrow \infty }\frac{R(t)}{t}=\mu^{-1} \text{,}  \label{1}
\end{equation}
where $\mu =ET_{1}\leq \infty.$ 

When $\mu <\infty $, Blackwell's theorem provides conditions under which%
\begin{equation}
\lim_{t\rightarrow \infty }\left( R(t+h)-R(t)\right) =\frac{h}{\mu }\text{,}
\label{2}
\end{equation}%
for all $h>0$. In the literature many papers are devoted to studying rates
of convergence in (1), (2).

In our paper we propose a general renewal process constructed by counting the number of successes associated with the generalized random walk $\{S_n \colon n \in \mathbb{N}_0\}$ with renewal function:
$$
R(t)=\sum_{n=0}^{\infty }P(S_{n}\leq t)=\sum_{n=0}^{\infty }F^{\diamond n}(t).
$$
We show that analogue of the Fredholm theorem holds for the renewal measure 
function associated with $\diamond$ and with the underlying unit step distribution $\nu$:
$$
m = \nu + \nu \diamond m.
$$
We prove an elementary renewal theorem for Kendall random walks in the case where the unit step d.f. is regularly varying. In particular, if $m(\alpha) = ET_1 < \infty$, then we arrive at
$$
\lim\limits_{t \to \infty} \frac{R(t)}{t^{\alpha
}} = \frac{2}{m(\alpha)}.
$$
We also prove a Blackwell theorem and a limit theorem for the renewal counting process in the Kendall random walk. For unit steps with finite $\alpha$-moment we have
$$
\lim_{t\rightarrow \infty} \left( \frac{R(t+h)}{(t+h)^{\alpha -1}} - \frac{R(t)}{t^{\alpha -1}} \right)= \frac{2h}{m(\alpha)}.
$$
All results on regularly varying cumulative distribution functions of
unit steps are related to its Williamson transform and its
asymptotic behavior.

\vspace{3mm}
Notation and organization of the paper is given as follows.
In Section 2 we give definitions
and properties of generalized convolutions. We recall the monotonicity property of generalized convolutions (for details see \cite{Poisson}) in the lack of memory property context. Next, we define renewal processes for generalized random walks. The principal result of Section 3 is Fredholms' theorem for monotonic renewal processes. We focus our investigations on the Kendall convolution case. We present a comprehensive result on the Williamson transform, which is the main mathematical tool in the Kendall convolution algebra. This transform corresponds to the characteristic function for classical convolution and has similar properties. The great advantage of the transform is that it is easy to invert for continuous probability measures. We find that the distribution of the constructed renewal process and corresponding moments. Most results are given in terms of the Williamson transform of the unit step and corresponding cumulative distribution function. Using regularly varying functions we prove Blackwell theorem and limit theorem for renewal
process constructed by the Kendall random walks with unit steps having regularly varying densities.

\vspace{3mm}
Throughout this paper, the distribution of the random
element $X$ is denoted by $\mathcal{L}(X)$. For a probability measure $\lambda$ and $a \in
\mathbb{R}_+$ the rescaling operator is given by $\mathbf{T}_a \lambda = \mathcal{L}(aX)$ if $\lambda = \mathcal{L}(X)$. For abbreviation the family of all probability measures on the Borel subsets of $\mathbb{R}_+$ is denoted by $\mathcal{P}_+$.

\vspace{3mm}

We consider the Kendall random walk $\{S_n \colon n \in \mathbb{N}_0\}$ defined in the following way:
\begin{defn} The stochastic process $\{S_n \colon n \in \mathbb{N}_0\}$ is a discrete time Kendall random walk  with parameter $\alpha>0$ and unit step distribution $\nu \in \mathcal{P}_+$ if there exist
\begin{namelist}{ll}
\item[\bf 1.] $(T_k)$ i.i.d. random variables with distribution $\nu$,
\item[\bf 2.] $(U_k)$ i.i.d. random variables with uniform distribution on $[0,1]$,
\item[\bf 3.] $(\theta_k)$ i.i.d. random variables with the symmetric Pareto distribution with the density $\pi_{2\alpha}\, (dy) = 2\alpha y^{-2\alpha-1} \pmb{1}_{[1,\infty)}(y)\, dy$,
\item[\bf 4.]  sequences $(T_k)$, $(U_k)$ and $(\theta_k)$ are independent,
\end{namelist}
such that
$$
S_0 = 1, \quad S_1 = T_1, \quad  S_{n+1} = M_{n+1}\left[ \mathbf{I}(U_n > \varrho_{n+1}) + \theta_{n+1} \mathbf{I}(U_n < \varrho_{n+1})\right],
$$
where $\theta_{n+1}$ and $M_{n+1}$ are independent and
$$
M_{n+1} = \max\{ S_n, T_{n+1}\}, \quad m_{n+1} = \min\{ S_n, T_{n+1}\}, \quad \varrho_{n+1} = \frac{m_{n+1}^{\alpha}}{M_{n+1}^{\alpha}}.
$$
\end{defn}

The process  can be also defined as a  random walk with respect to the Kendall generalized convolution as is done in the next section.

\section{Generalized convolutions and their properties}

\subsection{Generalized convolutions}
The following definition was introduced by Urbanik \cite{Urbanik64} influenced by the paper of Kingman \cite{King}. This was the beginning of the theory of generalized convolutions.
\begin{defn}
A generalized convolution is a binary, symmetric, associative and commutative operation $\diamond$ on $ \mathcal{P}_+$ with the following properties:
\begin{itemize}
\item[(i)] $\lambda \diamond \delta_0 = \lambda$ for all $\lambda \in \mathcal{P}_+$; \\
\item[(ii)] $(p\lambda_1 +(1-p)\lambda_2) \diamond \lambda = p(\lambda_1 \diamond \lambda) + (1-p)(\lambda_2 \diamond \lambda)$ for each $p \in [0,1]$ and $\lambda, \lambda_1, \lambda_2 \in \mathcal{P}_+$;\\
\item[(iii)] $\mathbf{T}_a(\lambda_1 \diamond \lambda_2) = (\mathbf{T}_a\lambda_1) \diamond (\mathbf{T}_a\lambda_2)$ for all $a \geq 0$ and $\lambda_1, \lambda_2 \in \mathcal{P}_+$;\\
\item[(iv)] if $\lambda_n \Longrightarrow \lambda$ and $\nu_n \Longrightarrow \nu$, then $(\lambda_n \diamond \nu_n) \Longrightarrow (\lambda \diamond \nu)$, where $\Longrightarrow $ denotes weak convergence;\\
\item[(v)] there exists a sequence of positive numbers $(c_n)$ such that $\mathbf{T}_{c_n} \delta_1^{\diamond n}$ converges weakly to a measure $\nu \neq \delta_0$ (here $\lambda^{\diamond n} = \lambda \diamond \lambda ... \diamond \lambda$ denotes the generalized convolution of $n$ identical measures $\lambda$).
\end{itemize}
\end{defn}
\begin{rem}
The pair $(\mathcal{P}_+, \diamond)$ is called a generalized convolution algebra. A continuous
mapping $h:\mathcal{P}_+\to \mathbb{R}$ is called a homomorphism of the algebra $(\mathcal{P}_+,\diamond)$ if
\begin{itemize}
\item[(i)] $h(a\lambda_1 + (1-a)\lambda_2) = ah(\lambda_1) + (1-a)h(\lambda_2)$ for every $a \in [0,1]$ and $\lambda_1, \lambda_2 \in \mathcal{P}_+$;
\item[(ii)] $h(\lambda_1 \diamond \lambda_2) = h(\lambda_1)h(\lambda_2)$ for every $\lambda_1, \lambda_2 \in \mathcal{P}_+$
\end{itemize}
Obviously, $h(\cdot) \equiv 0$ and $h(\cdot) \equiv 1$ are the trivial homomorphisms. A generalized convolution algebra $(\mathcal{P}_+, \diamond)$ is said to be regular if it admits existence of a non-trivial homomorphism.
\end{rem}
\begin{defn}
We say that a regular generalized convolution algebra $(\mathcal{P}_+, \diamond)$ admits a characteristic function if there exists a one-to-one
correspondence $\lambda \leftrightarrow \Phi_\lambda$ between probability measures $\lambda \in \mathcal{P}_+$ and real valued functions on $[0, \infty)$ such that for $\lambda, \nu, \lambda_n \in \mathcal{P}_+$
\begin{itemize}
\item[(i)] $\Phi_{p\lambda + q\nu} = p\Phi_\lambda + q\Phi_\nu$ for $p,q > 0$, $p+q=1$;\\
\item[(ii)] $\Phi_{\lambda \diamond \nu} =\Phi_\lambda \cdot \Phi_\nu$;\\
\item[(iii)] $\Phi_{\mathbf{T}_a\lambda}(t) = \Phi_\lambda(at)$; \\
\item[(iv)] the uniform convergence of $\Phi_{\lambda_n}$ on every bounded interval is equivalent to the weak convergence of $\lambda_n$.
\end{itemize}
The function $\Phi_\lambda$ is called the characteristic function of the probability
measure $\lambda$ in the algebra $(\mathcal{P}_+, \diamond)$ or the $\diamond$-generalized characteristic function of $\lambda$. 
\end{defn}

It is known that $\Phi_{\lambda}$ plays the same role as Fourier or Laplace transform for the classical convolution. It is important that in the regular generalized convolution algebra $(\mathcal{P}_+, \diamond)$ with the non-trivial homomorphism $h$ we have (up to a scale coefficient)
$$
\Phi_{\lambda} (t) = \int_0^{\infty} h(\delta_{xt}) \lambda(dx).
$$
The function $\Omega \colon [0,\infty) \rightarrow \mathbb{R}$ given by $\Omega(t) = h(\delta_t)$ is called the probability kernel for the generalized convolution $\diamond$. More information on homomorphism, generalized characteristic function and the probability kernel in a regular generalized convolution algebra one can find in \cite{Bin84}, \cite{MOU},  \cite{Urbanik64}.
\begin{Ex} \label{ex:alpha}
The stable $\alpha$-convolution, $\alpha >0$,  is defined for $a,b \geqslant 0$ by
$$
\delta_a \ast_{\alpha} \delta_b \stackrel{def}{=} \delta_c, \quad \hbox{ where } c^{\alpha} = a^{\alpha} + b^{\alpha}.
$$
In this case we have $h(\delta_x) = e^{-x^{\alpha}}$ and the corresponding generalized characteristic function 
$$
\Phi_{\lambda}(t) = \int_0^{\infty} e^{- t^{\alpha} x^{\alpha}} \lambda(dx).
$$
\end{Ex}
\begin{Ex} \label{ex:Kingman}
The Kingman convolution is easy to describe using independent random variables: If $\lambda_1 = \mathcal{L}(\theta_1)$, $\lambda_2 = \mathcal{L} (\theta_2)$ and the nonnegative variables $\theta_1, \theta_2$ are independent then
$$
\lambda_1 \otimes_n \lambda_2 = \mathcal{L} \left( \| \theta_1 U_1 + \theta_2 U_2 \|_2\right),
$$
where
$U_1, U_2$ are independent copies of the random vector $U$ with uniform distribution on the unit sphere in $\mathbb{R}^n$ and $\| \cdot \|_2$ denotes the Euclidean norm in $\mathbb{R}^n$. The corresponding generalized characteristic function is given by
$$
\Phi_{\lambda} (t) = \int_{\mathbb{R}} \Omega_n(tx) \lambda(dx),
$$
where $\Omega_n$ is the classical characteristic function (or, simply the Fourier transform) for the first coordinate of the random vector $U$.
\end{Ex}
\begin{Ex} \label{ex:Kendall}
The Kendall convolution $\vartriangle_{\alpha}$, $\alpha > 0$, is defined for $x\in [0,1]$ by
$$
\delta_x \vartriangle_{\alpha} \delta_1 = x^{\alpha} \delta_1 + (1-x^{\alpha}) \pi_{2\alpha},
$$
where $\pi_{\beta}$ denotes the Pareto distribution with density $\beta y^{-\beta-1} \mathbf{1}_{[1,\infty)}(y)$ (see \cite{Arnold2}). In this case we have $h(\delta_x) = (1-x^{\alpha})_+$, where $a_+ = a$ if $a>0$ and $a_+ = 0$ if $a \leqslant 0$. The corresponding generalized characteristic function is the Williamson transform (for more details see \cite{Williamson})
$$
\Phi_{\lambda}(t) = \int_0^{\infty} (1-x^{\alpha} t^{\alpha})_+ \lambda(dx).
$$
\end{Ex}
\begin{Ex} \label{ex:Kucharczak-Urbanik}
The Kucharczak-Urbanik convolution (\see e.g. \cite{JM_l-c}) is given by 
$$
\delta_x \diamond \delta_1 (ds) = (1 - x^{\alpha})^{n}\delta_1(ds) + \sum_{k=1}^n {n \choose k} x^{\alpha k} (1 - x^{\alpha})^{n-k} f_{k,n}(s) ds,
$$
where the densities $f_{k,n}$ are given by the following formulas
$$
f_{k,n}(s) = \alpha k {{n+k}\choose n} s^{-\alpha(n+1) - 1} \left( 1 - s^{-\alpha} \right)_+^{k-1}.
$$
The corresponding generalized characteristic function has the following expression
$$
\Phi_{\lambda}(t) = \int_0^{\infty} (1-x^{\alpha} t^{\alpha})^n_+ \lambda(dx).
$$
\end{Ex}
\begin{Ex} \label{ex:Civita}
For every $p \geqslant 2$ and properly chosen $c>0$  the function $h(\delta_t) = \varphi(t) = \varphi_{c,p}(t) = ( 1 - (c+1)t +ct^p)\mathbf{1}_{[0,1]}(t)$ is the kernel of a Kendall type (see \cite{JM_l-c}) generalized  convolution $\diamond$ defined for $x \in [0,1]$ by the formula:
$$
\delta_x \diamond \delta_1 = \varphi(x) \delta_1 + x^p \lambda_1 + (c+1)(x-x^p) \lambda_2,
$$
where $\lambda_1, \lambda_2$ are probability measures absolutely continuous with respect to the Lebesgue measure and independent of $x$. For example if $c = (p-1)^{-1}$ then
$$
\lambda_1(du) = \frac{2c}{u^3} \Bigl[ (c+1)(p+1) u^{1-p} + (c+1)(p-2) + cp (2p-1) u^{-2p-2} \Bigr] \!\!\mathbf{1}_{[1,\infty)}(u) du,
$$
and
$$
\lambda_2(du) = c \bigl[ 2(p-2) + (p+1)u^{-p+1} \bigr] u^{-3} \mathbf{1}_{[1,\infty)}(u) du.
$$
\end{Ex}
\begin{defn}
Let $\diamond$ be a generalized convolution on $\mathcal{P}_{+}$. The Markov process $(S_n)$ with $\mathcal{L}(S_1) = \nu$ and transition probabilities
$$
P_{k,n}(x,\cdot )=P(S_n\in \cdot\,|S_k=x):= \delta_x \diamond \nu^{\diamond (n-k)}(\cdot\,)
$$
we call a generalized random walk with respect to the convolution $\diamond$ with the step distribution $\nu$.
\end{defn}

\subsection{Monotonicity property and lack of memory property}

For the classical renewal process the corresponding process of change times $T_1 + \dots T_n$, $n\in \mathbb{N}$, is monotonically increasing. It turns out  that in the case of generalized convolutions the sequence $S_n$ does not have to be monotonically increasing. For this reason the monotonicity property for generalized convolutions on $\mathcal{P}_{+}$ was defined in \cite{Poisson} in the following way:
\begin{defn}
A generalized convolution $\diamond$ is monotonic if for all $x,y \geqslant 0$, $x\vee y =\max\{x,y\}$
$$
\delta_x \diamond \delta_y \left( [x\vee y, \infty) \right) = 1.
$$
\end{defn}
Notice that classical convolution, stable $\ast_{\alpha}$ convolution and the Kendall $\vartriangle_{\alpha}$ convolution are monotonic. The Kingman convolution $\otimes_n$ is not monotonic since the measure $\delta_x \otimes_n \delta_y$ is supported in $[|x-y|, x+y]$.

For the proofs of the results below see \cite{Poisson}.
By $F_n$ we denote the cumulative distribution function of the variable $S_n$.

\begin{lem}\label{lem:2}
If the generalized convolution $\diamond$ is monotonic then for every $k_1< k_2< \dots< k_{n+1}$
\begin{eqnarray*}
{\rm a)}  & \mathbf{P}\left\{ S_{k_{n+1}} > x, S_{k_n} > x , \dots, S_{k_1}> x\right\} = 1-F_{k_1}(x), \\
{\rm b)} & \mathbf{P}\left\{ S_{k_{n+1}} > x, S_{k_n} \leqslant x, \dots, S_{k_1} \leqslant x \right\} = \mathbf{P}\left\{ S_{k_{n+1}} > x, S_{k_n} \leqslant x \right\} \\
 & = F_{k_n}(x) - F_{k_{n+1}}(x).
\end{eqnarray*}
\end{lem}

The classical construction of the Poisson process $\{ N(t) \colon t \geqslant 0\}$ is based on a sequence $(T_k)$ of i.i.d. random variables with the exponential distribution $\Gamma(1,a)$ which have the lack of memory property. It turns out that this property strongly depends on the considered convolution: 

\begin{defn}
A probability distribution $\nu$ with the cumulative distribution function $F$ has the lack of the memory property with respect to generalized convolution $\diamond$ if
$$
\mathbf{P}\bigl\{ X > x \diamond y  \big| X > x \bigr\} =  1 - F(y), \quad x,y \in \mathbb{R}_{+},
$$
where $X$ with distribution function $F$ is independent of the random variable $x \diamond y$ having distribution $\delta_x \diamond \delta_y$.
\end{defn}

\vspace{2mm}

\begin{prop}\label{monoton}
The distribution function $F\not\equiv \mathbf{1}_{[0,\infty)}$ has the lack of memory property with respect to the monotonic generalized convolution $\diamond$ if and only if the algebra $(\mathcal{P}_{+}, \diamond)$ is regular with homomorphism $h(\delta_t)$, $t \geqslant 0$, which is monotonically decreasing as a function of $t$ and $F(t) = 1 - h(\delta_{c^{-1}t})$ for some $c>0$.
\end{prop}

\vspace{2mm}

\noindent
\begin{rem}
{\rm  It follows from Proposition \ref{monoton} that the regular generalized convolution $\diamond$ admits the existence of a measure  with the lack of memory property with respect to $\diamond$ if and only if the  homomorphism $h$ is the tail of some distribution function. }
\end{rem}

\noindent
\begin{rem}
{\rm It is evident now that the stable convolution and the Kendall convolution admit the existence of a distribution with the lack of memory property since their homomorphisms $h(\delta_t) = e^{-t^{\alpha}}$ and $h(\delta_t) = (1-t^{\alpha})_{+}$ respectively are the tails of some distribution functions.}
\end{rem}

\section{Renewal processes with respect to generalized convolution}
In this section we focus on the renewal theory and survey some of the most important developments in probability theory. The classical theory of renewal processes has been extensively studied in \cite{Asmussen}, \cite{MitovOmey} and \cite{Rol99}. 
Some generalizations of renewal processes, called Wold processes, one can find in \cite{Daley}. In this paper we consider renewal processes for Markov chains with transition probabilities defined by generalized convolutions (see \cite{BJMR}) in a way they can be treated as special cases of the Wold processes.
\begin{defn}
Let $\diamond$ be a generalized convolution on $\mathcal{P}_+$, $\{T_n \colon n \in \mathbb{N} \}$ be a sequence of i.i.d. random variables with distribution $\nu$ and $\{S_n \colon n \in \mathbb{N}_0\}$ be the corresponding generalized random walk. Then the process $\{ N(t) \colon t \geqslant 0\}$ defined by
$$
N(t) = \left\{ \begin{array}{l}
 \inf\{ n \colon S_{n+1} > t \} \\[2mm]
 \infty \; \hbox{ if such $n$ does not exists }
 \end{array} \right.
$$
is called the $\diamond$-renewal process. If $\nu$ has the lack of memory property then the $\diamond$-renewal process is called the $\diamond$-Poisson process.
\end{defn}

Notice that for every generalized convolution $\diamond$ we have
$$
\mathbf{P}\{ N(t) = 0 \} = \nu (t,\infty)= 1 - F(t).
$$
\subsection{Fredholm theorem  for $\diamond$-renewal process}
\begin{defn}
Let $\{ N(t) \colon t \geqslant 0\}$ be a $\diamond$-renewal process for generalized convolution $\diamond$. The renewal function is defined by
$$
R(t) \stackrel{def}{=} \mathbf{E}N(t).
$$
\end{defn}
\begin{lem}
Let $\{N(t)\colon t \geqslant 0\}$ be a $\diamond$-renewal process with unit step distribution $\nu$. If the generalized convolution $\diamond$ is monotonic then  
$$
\mathbf{P}\{ N(t) = n \} = F_{n}(t) - F_{n+1}(t),
$$
where $F_n$ is the cumulative distribution function of $S_n$. Moreover
$$
R(t) = \sum_{n=1}^{\infty} F_n(t), \quad \mathbf{E}N^2(t) = 2 \sum_{n=1}^{\infty} n F_n(t) - \sum_{n=1}^{\infty} F_n(t).
$$
\end{lem}

In order to prove an analogue of the Fredholm renewal equation we introduce the following generalized convolution transform:
\begin{defn}
 Let $h(\delta_t)$ be the probability kernel for the algebra $(\mathcal{P_+}, \diamond)$. The generalized convolution $\diamond$ transform for the function $f: \mathbb{R}_+ \to \mathbb{R}_+ $ is defined by:
$$
\widehat{f}(t) = \int_0^\infty h(\delta_{ts})f(s)ds.
$$
\end{defn}
Notice that if  $T \sim \nu$ has density function $f$, then the generalized characteristic function $\Phi_\nu$ can be expressed by $\widehat{f}$:
$$
\Phi_\nu(t) \stackrel{def}{=} \int_0^\infty h(\delta_{ts}) \nu(ds) = \int_0^\infty h(\delta_{ts})f(s)ds = \widehat{f}(t).
$$

\begin{thm}{\bf Fredholm theorem for generalized convolution.}
Let $R(t)$ be the renewal function for the $\diamond$-renewal process  $\{N(t)\colon t \geqslant 0\}$ with unit step distribution $\nu$, where the generalized convolution $\diamond$ is monotonic. Then $R(t)$ is the cumulative distribution function of a non-negative measure $m$ and:
$$
m = \nu + \nu \diamond m.
$$
\end{thm}

\noindent
{\bf Proof.}
From the definition of the $\diamond$-renewal process we have:
$$
R(t) = \mathbf{E}N(t) = \sum_{n=1}^{\infty} F_n(t) = \sum_{n=1}^{\infty} F_{\nu^{\diamond n}}.
$$
Since $R$ is a sum of cumulative distribution functions it is also a cumulative distribution function of a positive measure $m$. 
Using generalized characteristic functions, we calculate:
\begin{eqnarray*}
\Phi_m (t) & = & \int_0^\infty h(\delta_{ts}) dR(s)  =  \sum_{n=1}^\infty \int_0^\infty h(\delta_{ts})dF_n(s) \\
& = & \sum_{n=1}^\infty \Bigl\{ \int_0^\infty h(\delta_{ts})\nu(ds) \Bigr\}^n \\
& = & \int_0^\infty h(\delta_{ts}) \nu(ds) + \int_0^\infty h(\delta_{ts}) \nu(ds) \int_0^\infty h(\delta_{ts})dR(s)\\
& = & \Phi_\nu(t) + \Phi_\nu(t) \Phi_m(t).
\end{eqnarray*}
This finally implies that $m = \nu + \nu \diamond m$. Moreover we have $\Phi_m(t) = {{\Phi_{\nu}(t)}/{(1-\Phi_{\nu}(t))}}$. \qed

\vspace{2mm}

Notice that  the formula $m = \nu + \nu \ast m$ is the simplified notation for the classical Fredholm result (with the classical convolution $\ast$): 
$$
m(t) = F(t) + \int_0^t m(t-x) dF(x).
$$

\vspace{2mm}

\subsection{Renewal process with respect to the Kendall convolution}
Henceforth we focus our investigations on the Kendall convolution case, which was defined in Example \ref{ex:Kendall} in the following way:
$$
\delta_x \vartriangle_\alpha \delta_1 = (1-x^{\alpha}) \delta_1 + x^{\alpha} \pi_{2\alpha},
$$
for $0<x\leq 1$ and $\alpha>0$, where $\pi_{2\alpha}(dx) = \frac{2\alpha}{x^{2\alpha-1}}\mathbf{1}_{[1,\infty)} dx$.

The  Kendall convolution probability kernel is given by
$$
h(x,y,t) := \delta_x \vartriangle_{\alpha} \delta_y (0,t) = \left( 1 - \frac{x^{\alpha}y^{\alpha}}{t^{2\alpha}}\right) \mathbf{1}_{\{ x<t,\, y<t\}}.
$$
for \,$x,y >0$.

The main tool, which we use in the Kendall convolution algebra, is the Williamson transform (see \cite{KendallWalk},\cite{factor}), the operation $\nu \rightarrow \widehat{\nu}$ given by
$$
\widehat{\nu}(t) = \int_{\mathbb{R}_+} \left( 1 - (xt)^{\alpha} \right)_+ \nu(dx), \quad \nu \in \mathcal{P}_+,
$$
where $a_+ = a$ if $a\geqslant 0$ and $a_+=0$ otherwise. The function $\widehat{\nu}$ is the generalized characteristic function in the Kendall convolution algebra. For convenience we use the following notation:
$$
G(t) \stackrel{def}{=} \widehat{\nu}(1/t) = \int_0^{\infty} \Bigl( 1 - \frac{x^{\alpha}}{t^{\alpha}} \Bigr)_{+} \nu(dx), \quad \Psi(t) = \left( 1 - t^{\alpha} \right)_+.
$$
The Williamson transform has the following important  property:
\begin{prop}\label{prop:2}
Let $\nu_1, \nu_2 \in \mathcal{P}_+$ be probability measures with Williamson transforms $\widehat{\nu_1}, \widehat{\nu_2}$. Then
$$
\int_{\mathbb{R}_+} \Psi(xt) \bigl(\nu_1 \vartriangle_{\alpha} \nu_2 \bigr)(dx) = \widehat{\nu_1}(t) \widehat{\nu_2}(t).
$$
\end{prop}
The next lemma describes the inverse of the Williamson transform.  
\begin{lem}\label{lem:3}
Let $G(t)= \widehat{\nu}(t^{-1})$ and $H(t) = \int_0^t x^{\alpha} \nu(dx)$ for the probability measure $\nu \in \mathcal{P}_+$ with cumulative distribution function  $F$, $\nu(\{0\}) = 0$. Then
$$
F(t) =  G(t)  + \frac{t}{\alpha} G'(t) = G(t)  + t^{-\alpha} H(t)
$$
except, possibly,  for countable many points $t \in \mathbb{R}$.
\end{lem}
\noindent
{\bf Proof.} It is enough to notice that
\begin{eqnarray*}
G(t) + \frac{t}{\alpha}G'(t) &=& \int_0^t(1-x^\alpha t^{-\alpha})\nu(dx) + \frac{t}{\alpha}\int_0^t \alpha x^\alpha t^{-\alpha-1}\nu(dx)  \\
&=& \int_0^t \nu(dx) = F(t). 
\end{eqnarray*}
The equality $G'(t) = \alpha t^{-\alpha - 1} H(t)$ is trivial. \qed

\vspace{2mm}

Using Lemma \ref{lem:3} we see that $G'(t) = \frac{\alpha}{t} \left(F(t)-G(t)\right)$. The Williamson transform corresponding to $F_n$ is given by $G_n(t) = G(t)^n$. Hence we have the following:

\begin{lem}\label{lem:4}
Let $\diamond$ be the Kendall generalized convolution, $\{N(t) \colon t \geqslant 0 \}$ be the $\diamond$-renewal process with unit step distribution $\nu$ and let $F_n$ be the cumulative distribution function of $\nu^{\diamond n}$. Then
\begin{eqnarray*}
F_n(t) & = & G(t)^{n-1} \bigl[n\left(F(t) - G(t)\right) + G(t) \bigr] \\
& = & G(t)^{n-1} \bigl[n t^{-\alpha} H(t) + G(t) \bigr] .
\end{eqnarray*}
Consequently
\begin{eqnarray*}
\mathbf{P} \left\{ N(t) = 0 \right\}\!\! & = & \!\! \overline{F}(t), \\ 
\mathbf{P} \left\{ N(t) = n \right\}\!\! & = & \!\!G(t)^{n-1} \bigl[ n (F(t) - G(t))\overline{G}(t) + G(t) \overline{F}(t)) \bigr],
\end{eqnarray*}
where
$$
\overline{F}(t) = 1-F(t), \quad \overline{G}(t) = 1-G(t).
$$
\end{lem}

Now we are able to calculate the renewal function for the renewal process with respect to the Kendall generalized convolution:
\begin{thm}
Let $\{ N(t) \colon t \geqslant 0\}$ be a $\diamond$-renewal process for the Kendall convolution $\diamond = \triangle_{\alpha}$, $\alpha >0$. The renewal function is given by
$$
R(t) \stackrel{def}{=} \mathbf{E} N(t) = \frac{G(t)}{\overline{G}(t)} + \frac{H(t)}{t^{\alpha}\overline{G}(t)^2}.
$$
Moreover
$$
\mathbf{E}N^2(t) = \frac{t^{-\alpha} H(t) \left(1+3G(t)\right) +  G(t) \left(1-G^2(t)\right)}{\overline{G}(t)^3}
$$
and
$$
Var\left( N(t) \right) = \frac{H(t) \left(1 + G(t)\right)}{t^{\alpha} \overline{G}(t)^3} + \frac{G(t)}{ \overline{G}(t)^2} - \frac{H^2(t) }{t^{2\alpha} \overline{G}(t)^4}.
$$
\end{thm}
\noindent
{\bf Proof.} In order to prove the formulas we shall use the following properties of the geometrical series:
$$
\sum\limits_{n=1}^{\infty} nq^{n-1} =\frac{1}{(1-q)^2}, \quad
\sum\limits_{n=1}^{\infty} n^2q^n = \frac{q^2+q}{(1-q)^3}.
$$
By the second expression for $F_n$ given in Lemma \ref{lem:4} we have:
\begin{eqnarray*}
R(t) & = & \sum_{n=1}^{\infty} F_n(t) = \sum_{n=1}^{\infty} G(t)^{n} + t^{-\alpha}H(t) \sum_{n=1}^{\infty} n G(t)^{n-1}  \\
& = & \frac{G(t)}{\overline{G}(t)} + \frac{H(t)}{t^{\alpha}\overline{G}(t)^2},
\end{eqnarray*}
and
\begin{eqnarray*}
\sum_{n=1}^{\infty}nF_n(t) & = & \sum_{n=1}^{\infty}n\Bigl[G^n(t) + n\left(F(t)-G(t)\right)G^{n-1}(t)\Bigr] \\
& = & \frac{G(t)\left(1-G(t)\right) + \left(F(t)-G(t)\right)(G(t)+1)}{\overline{G}(t)^3},
\end{eqnarray*}
which yields
\begin{eqnarray*}
\mathbf{E}N^2(t) & = & 2 \sum_{n=1}^{\infty} n F_n(t) - \sum_{n=1}^{\infty} F_n(t) \\
& = & \frac{t^{-\alpha} H(t) \left(1+3G(t)\right) + G(t) \left(1-G^2(t)\right)}{\overline{G}(t)^3}.
\end{eqnarray*}
Now using expression for $R^2(t)$ we arrive at:
$$
Var(N(t))  =   \frac{\left(F(t)-G(t)\right)\left(1 + G(t)\right)}{\overline{G}(t)^3} + \frac{G(t)}{ \overline{G}(t)^2} - \frac{\left(F(t)-G(t)\right)^2 }{ \overline{G}(t)^4}. \eqno{\Box}
$$
\vspace{2mm}

\begin{Ex}\label{ex:delta1}
Let $\nu=\delta_1$. Then $m(\alpha) = {\mathbf E} S_1^{\alpha} = 1$, 
$$
G(x) = (1- x^{-\alpha})_+, \quad H(x) = \mathbf{1}_{[1,\infty)}(x),
$$
and
$$
R(t) = (2t^{\alpha} -1)\mathbf{1}_{[1,\infty)}(t).
$$
\end{Ex}

\begin{Ex}\label{ex:uniform}
Let $\nu=U(0,1)$ be the uniform distribution on the interval $(0,1)$.
Then $m(\alpha) = {\mathbf E} S_1^{\alpha} = (\alpha+1)^{-1}$ and the Williamson transform is the following:
$$
G(x) = (x \wedge 1) - \frac{(x \wedge 1)^{\alpha+1}}{(\alpha+1)x^{\alpha}},\quad H(x) = \frac{(x \wedge 1)^{\alpha+1}}{(\alpha+1)},
$$
$$
R(t) = \left\{ \begin{array}{lcl}
 t \, \frac{(\alpha+1)^2-\alpha^2t}{(\alpha+1-\alpha t)^2} & \hbox{ if } & t \leq 1; \\
 2(\alpha+1)t^\alpha-1 & \hbox{ if } & t > 1.
 \end{array} \right.
$$
\end{Ex}

\begin{Ex}\label{ex:pareto}
Let $\nu = \pi_{2\alpha}$, where $\alpha>0$ and $\pi_{p}$ be the Pareto distribution with the density
$$
\pi_{p}(dx)=\frac{p}{x^{p+1}} \pmb{1}_{[1,\infty)}(x)dx.
$$
Then $m(\alpha) = {\mathbf E} S_1^{\alpha} = 2$,
$$
G(x) = (1-x^{-\alpha})_+^2, \quad H(x)= 2(1-x^{-\alpha})_+
$$
and
$$
R(t) = \frac{(t^{\alpha}- 1)(1 - 3t^{\alpha} + 4t^{2\alpha})}{(2t^{\alpha} -1)^2}\, {\mathbf 1}_{(1,\infty)}(t).
$$
\end{Ex}

\begin{Ex}\label{ex:lackmem}
For the Kendall convolution (\cite{MisJas2}) the distribution having the tail:
$$
\nu(x,\infty) = (1-x^{\alpha})_+
$$
has the lack of memory property.
Then $m(\alpha) = {\mathbf E} S_1^{\alpha} = \frac{1}{2}$,
$$
G(x) = \frac{x^{\alpha}}{2}\pmb{1}_{[0,1)}(x) + \left(1-\frac{1}{2x^{\alpha}}\right) \pmb{1}_{[1,\infty)}(x), \quad
H(x) = \frac{(x \wedge 1)^{2\alpha}}{2},
$$
$$
R(t)= \left\{ \begin{array}{lcl}
 (4t^{-\alpha}-1)(2t^{-\alpha}-1)^{-2} & \hbox{ if } & t \leq 1; \\[1mm]
 4t^{\alpha}-1 & \hbox{ if } & t > 1.
 \end{array} \right.
$$
\end{Ex}

\begin{Ex}\label{ex:stable}
As a step distribution we consider the $\alpha$-stable distribution at the Kendall convolution algebra with the generalized characteristic function $\Phi(t) = e^{-t^{\alpha}}$ on $[0,\infty)$ and the following density function:
$$
\nu(dx) = \alpha x^{-(2\alpha +1)} exp\{-x^{-\alpha}\} \mathbf{1}_{[0,\infty)}(x)dx.
$$ 
Then $m(\alpha) = {\mathbf E} S_1^{\alpha} = 1$,
$G(x) =  e^{-x^{-\alpha}} = H(x)$.
Consequently the renewal function is given by:
$$
R(t) = \bigl(e^{t^{-\alpha}}-1+t^{-\alpha}e^{t^{-\alpha}}\bigr)\bigl(e^{t^{-\alpha}}-1 \bigr)^{-2}.
$$
\end{Ex}

\vspace{2mm}

\subsection{Blackwell theorem for the Kendall renewal process}

It is interesting to study the asymptotic behaviour of $R(t)$ and that of the difference $R(t+h)-R(t)$. The result regarding the limit behavior of this difference is given in one of the most important theorems in the classical renewal theory - the Blackwell theorem. It states that $R(t+h)-R(t)$ converges to a constant as $t \to \infty$. For details see \cite{MitovOmey} and \cite{Omey}.

In the renewal theory for Kendall random walks this difference doesn't converge unless $\alpha=1$ so the analogue of the Blackwell theorem has a more complicated form.

In the following we use the idea of regularly varying functions at infinity (for a survey of regularly varying functions and their applications we refer to \cite{BGT}, \cite{Seneta}).

\begin{defn}
A positive and measurable function $f$ is regularly varying at infinity and
with index $\beta $ (notation $f\in RV_{\beta }$) if it satisfies
\[
\lim_{t\rightarrow \infty }\frac{f(tx)}{f(t)}=x^{\beta },\forall x>0\text{.}
\]
\end{defn}

Examples include functions such as $f(x)=x^{\beta }\log x$, $
f(x)=(1+x^{\alpha })^{\beta }$, $f(x)=\log \log x$, and so on. Notice that if $f(t) \to c$ for $t\to \infty$, then $f(tx)/f(t) \to 1$ for $t\to \infty$ and then $f\in RV_0$. Consequently in our case if $\mathbf{E}X^\alpha < \infty$ then $H\in RV_0$.

\vspace{2mm}

The next result is known as Karamata theorem (\cite{BGT}, Sect. 1.6,
p. 26).

\begin{thm}
Suppose that $f\in RV_{\beta }$ and that $f$ is bounded on bounded intervals.

(i) If $\beta \geq -1$, then $xf(x)/\int_{0}^{x}f(t)dt\rightarrow \beta +1$;

(ii) If $\beta <-1$, then $xf(x)/\int_{x}^{\infty }f(t)dt\rightarrow -(\beta
+1)$.
\end{thm}

\begin{lem}\label{lem:5}
Suppose that $H\in RV_{\theta }$ where $0\leq \theta <\alpha $ and  $\alpha W(x) =  H(x) +  x^{\alpha} \overline{F}(x)$. Then

(i) $x^{\alpha }\overline{F}(x)/H(x)\rightarrow \theta /(\alpha -\theta )$, for $x \rightarrow \infty$;

(ii) $W(x)/H(x)\rightarrow \alpha/(\alpha -\theta )$, for $x \rightarrow \infty$;

(iii) $x^{\alpha }\overline{G}(x)/H(x)\rightarrow \alpha /(\alpha -\theta )$ , for $x \rightarrow \infty$.
\end{lem}

\noindent{\bf Proof.} 
(i) If $H\in RV_{\theta }$ where $0\leq \theta <\alpha $, then $x^{-\alpha
-1}H(x)\in RV_{\theta -\alpha -1}$ and Karamata's theorem gives that 
\[
\frac{x^{-\alpha }H(x)}{\int_{x}^{\infty }y^{-\alpha -1}H(y)dy}\rightarrow
\alpha -\theta
\]
for $x \rightarrow \infty$.
On the other hand, we have
\begin{eqnarray*}
\lefteqn{\int_{x}^{\infty }y^{-\alpha -1}H(y)dy = \int_{y=x}^{\infty }y^{-\alpha
-1}\int_{z=0}^{y}z^{\alpha }dF(z)dy } \\
&=&\int_{y=x}^{\infty }y^{-\alpha -1}\int_{z=0}^{x}z^{\alpha
}dF(z)dy+\int_{y=x}^{\infty }y^{-\alpha -1}\int_{z=x}^{y}z^{\alpha }dF(z)dy
\\
&=&\frac{H(x)x^{-\alpha }}{\alpha }+\int_{z=x}^{\infty }\int_{y=z}^{\infty
}y^{-\alpha -1}z^{\alpha }dF(z)dy \\
&=&\frac{H(x)x^{-\alpha }}{\alpha }+\frac{1}{\alpha }\int_{z=x}^{\infty
}dF(z)= \frac{H(x)x^{-\alpha }}{\alpha }+\frac{1}{\alpha }\overline{F}(x)\text{.}
\end{eqnarray*}
We find that
\[
\frac{x^{\alpha }\overline{F}(x)}{H(x)}\rightarrow 
\frac{\theta }{\alpha -\theta }\text{,} \quad x \rightarrow  \infty.
\]
To see ($ii$) note that by these calculations we have
$$
H(x)=\int_{0}^{x}y^{\alpha }dF(y)=\alpha \int_{0}^{x}y^{\alpha -1}
\overline{F}(y)dy-x^{\alpha }\overline{F}(x) = \alpha W(x) - x^{\alpha }\overline{F}(x).
$$
Consequently $\alpha W(x)/H(x)\rightarrow 1+\theta /(\alpha -\theta )=\alpha
/(\alpha -\theta )$ and the result follows. The property ($iii$) easily follows from ($i$) since 
$\overline{G}(x)=\overline{F}(x)+x^{-\alpha }H(x)$
\qed

The following result is the elementary renewal theorem.
\begin{cor}
Suppose that $H\in RV_{\theta }$ where $0\leq \theta <\alpha $. Then
\[
x^{-\alpha }R(x)H(x)\rightarrow \alpha ^{-2}(\alpha -\theta )(2\alpha
-\theta )
\]
for $x \rightarrow \infty$.
\end{cor}

\noindent
{\bf Proof.} 
We have
\[
x^{-\alpha }R(x)H(x)=G(x)H(x)/x^{\alpha }\overline{G}(x)+H^{2}(x)/x^{2\alpha
}\overline{G}^{2}(x)
\]
and from Lemma \ref{lem:5} it follows that
\[
x^{-\alpha }R(x)H(x)\rightarrow \frac{\alpha -\theta }{\alpha }+\left(\frac{
\alpha -\theta }{\alpha }\right)^{2}=\frac{(\alpha -\theta )(2\alpha -\theta )}{
\alpha ^{2}}\text{.}
\]
This proves the result.
\qed

\begin{cor}
If $H(\infty )=m(\alpha )<\infty $, then $H(x)\in RV_{0}$ and for $x\rightarrow \infty$ we have the following 
$x^{\alpha }\overline{F}(x)\rightarrow 0$, $
x^{\alpha }\overline{G}(x)\rightarrow m(\alpha )$ and $x^{-\alpha
}R(x)\rightarrow 2/m(\alpha )$.
\end{cor}

\vspace{2mm}

\textbf{Remarks.}
Suppose that $H\in RV_{\theta }$ where $0\leq \theta <\alpha $.

1) Since $\theta <\alpha $, we have $x^{-\alpha }H(x)\rightarrow 0$ and $
\overline{G}(x)\rightarrow 0$.

2) If $\theta >0$,then $\overline{F}(x)\in RV_{\theta -\alpha }$.

\bigskip
Our next result is a Blackwell type theorem. We assume here that $F$ has a density 
which we denote by $f(x)$. Using $H(x)=\int_{0}^{x}y^{\alpha }dF(y)$, we have that 
$H^{\prime }(x)=x^{\alpha}f(x)$.

\begin{thm}
Suppose the cumulative distribution function $F$ is differentiable with the density $f$ 
and let $\lim\limits_{t\to\infty} t H^{\prime }(t)/H(t)=\theta $, where $0\leq \theta
<\alpha $. Then we have
$$
\quad \lim_{t\rightarrow \infty} \frac{H(t)R^{\prime }(t)}{ t^{\alpha-1 }}= \left(\frac{\alpha
-\theta }{\alpha }\right)^{2}(2\alpha -\theta )=:c; \leqno{(i)}
$$
$$
\forall \, h>0 \quad \lim_{t\rightarrow \infty} \frac{H(t) (R(t+h)-R(t))}{t^{\alpha -1 } } = ch. \leqno{(ii)}
$$
\end{thm}

\noindent
{\bf Proof.} 
(i) We have
$$
R^{\prime }(t) = \frac{f(t)}{\overline{G}^{2}(t)} + \frac{2\alpha H^{2}(t)}{t^{2\alpha +1} 
\overline{G}^{3}(t)}  = \frac{H^{\prime }(t) }{t^{\alpha}\overline{G}^{2}(t)} + \frac{2\alpha
H^{2}(t)}{t^{2\alpha +1} \overline{G}^3 (t) }. 
$$
Using Lemma \ref{lem:5} we obtain
\begin{eqnarray*}
\lefteqn{t^{-\alpha +1}H(t)R^{\prime }(t) = \frac{tH^{\prime }(t)}{H(t)}
\frac{H^{2}(t)}{t^{2\alpha}\overline{G}^{2}(t)} + 2\alpha \frac{H^{3}(t)}{t^{3\alpha } 
\overline{G}^{3} (t)} } \\
&& \longrightarrow \theta \left(\frac{\alpha -\theta }{\alpha }\right)^{2}+2\alpha \left(\frac{
\alpha -\theta }{\alpha }\right)^{3} = \left(\frac{\alpha -\theta }{\alpha }\right)^{2}(2\alpha
-\theta ),
\end{eqnarray*}
when $t\rightarrow \infty$. The result follows.

To prove ($ii$), note that
\[
R(t+h)-R(t)=\int_{0}^{h}R^{\prime }(t+z)dz\text{.}
\]
Since $R^{\prime }(t)\sim ct^{\alpha -1}/H(t)\in RV_{\alpha -\beta -1}$ when $t\rightarrow \infty$, we
have $R(t+h)-R(t)\sim R^{\prime }(t)h$ when $t\rightarrow \infty$, and the result follows.
\qed

\vspace{2mm}

\begin{cor}
If $m(\alpha )<\infty $ and $xH^{\prime }(x)\rightarrow 0$, we have $
x^{1-\alpha }(R(x+h)-R(x))\rightarrow 2\alpha h/m(\alpha )$.
\end{cor}

\vspace{2mm}

In the following theorem we present another version of the Blackwell theorem for generalized Kendall convolution:

\begin{thm}
Under the conditions of Corollary 3, we have
$$
\lim_{t\rightarrow \infty} \left( \frac{R(t+h)}{(t+h)^{\alpha -1}} - \frac{R(t)}{t^{\alpha -1}} \right)= \frac{2h}{m(\alpha)}.
$$
\end{thm}

\noindent
{\bf Proof.}  
We consider $R(t) t^{1-\alpha}$. Clearly we have
$$
\bigl(R(x)x^{1-\alpha }\bigr)' = R'(x)x^{1-\alpha }+(1-\alpha )R(x)x^{-\alpha }.
$$
From Theorem 4 (with $\theta =0$), we have $R^{\prime }(x)x^{1-\alpha
}\rightarrow 2\alpha /m(\alpha )$ and in\ Corollary 2 we proved that $
R(x)x^{-\alpha }\rightarrow 2/m(\alpha )$.
We find that for $x \rightarrow \infty$
\[
\left(\frac{R(x)}{x^{\alpha -1}}\right)' \rightarrow \frac{2\alpha }{m(\alpha )} + \frac{
2(1-\alpha )}{m(\alpha )}=\frac{2}{m(\alpha )}.
\]
Using the well known Lagrange theorem on the mean value for the differentiable function $f$: 
$$
f(x+h) - f(x) = h f'(x + \varrho h) \quad \hbox{ for some } \quad \varrho \in [0,1]
$$ 
we conclude that for some $\varrho = \varrho(x) \in [0,1]$ and $x \rightarrow \infty$ 
\[
\frac{R(x+h)}{(x + h)^{\alpha-1}}-\frac{R(x)}{x ^{\alpha-1}} = h\, \left(\frac{R(x + \varrho h )}{(x + \varrho h) ^{\alpha-1} }\right)' \longrightarrow \frac{2h }{m(\alpha )}.  \eqno{\Box}
\]

\vspace{2mm}

\subsection{A limit theorem for the renewal process}

Now we consider the limiting behavior of $N(t)$.

\begin{thm}
(i) If $H\in RV_{\theta }$ where $0\leq \theta <\alpha $, then $\overline{G}(t)N(t)\Longrightarrow Z$, 
where the distribution of $Z$ is a convex linear combination of  an exponential and a gamma distribution:
$$
\mathcal{L}(Z) = \alpha^{-1} \theta \, \Gamma(1,1) +  (1- \alpha^{-1} \theta ) \,\Gamma(2,1),
$$
where $\Gamma(p,b)$ denotes the measure with the density $\frac{b^p}{\Gamma(p)} x^{p-1} e^{-bx} {\mathbf 1}_{[0,\infty)}(x)$. 
(ii) If $\theta =0$, then $\overline{G}(t)N(t)\Longrightarrow Z$ and $\mathcal{L}(Z) = \Gamma(2,1)$.
\end{thm}

\noindent{\bf Proof.}
Using $F_k(t) = G^{k-1}(t) \bigl(k(F(t)-G(t))+G(t)\bigr)$ and $\mathbf{P}\bigl(N(t)=k\bigr) = F_k(t) - F_{k+1}(t)$ we calculate
the Fourier transform of $N(t)$:
\begin{eqnarray*}
\lefteqn{\Phi_{N(t)}(u)
= \sum_{k=1}^\infty e^{iuk}\bigl(F_k(t)-F_{k+1}(t) \bigr) + \bigl(1-F(t) \bigr)} \\
& = & e^{iu}F(t) + \left(e^{iu} - 1 \right)\sum_{k=2}^\infty e^{iu(k-1)}F_k(t) + 1 - F(t) \\
& = & 1 + (e^{iu} - 1)F(t) + (e^{iu} - 1)\bigl(F(t)-G(t) \bigr)\sum_{k=1}^\infty k\bigl(e^{iu}G(t) \bigr)^{k-1}  \\
& - & (e^{iu}-1) \bigl(F(t)-G(t) \bigr) + (e^{iu}-1)e^{-iu}\sum_{k=1}^\infty \bigl(e^{iu}G(t) \bigr)^k - (e^{iu}-1)G(t) \\
& = & 1 + (e^{iu} - 1)\left[(F(t)-G(t))\sum_{k=1}^\infty k(e^{iu}G(t))^{k-1} + e^{-iu} \sum_{k=1}^\infty (e^{iu}G(t))^k \right] \\
& = & 1 + (e^{iu} - 1) \left[\bigl(F(t)-G(t) \bigr)\frac{1}{\bigl(1-e^{iu}G(t) \bigr)^2} + e^{-iu} \frac{e^{iu}G(t)}{1-e^{iu}G(t)} \right] \\
& = & 1 + (e^{iu} - 1)\frac{F(t) - e^{iu}G^2(t)}{(1-e^{iu}G(t))^2}.
\end{eqnarray*}
Consequently
\[
\Phi _{N(t)}(u)=1+(e^{iu}-1)\frac{F(t)-e^{iu}G^{2}(t)}{(1-e^{iu}G(t))^{2}}
\]
It follows that
\[
\Phi _{N(t)}(u)=1+(e^{iu}-1)\frac{F(t)-G^{2}(t)-(e^{iu}-1)G^{2}(t)}{(
\overline{G}(t)-(e^{iu}-1)G(t))^{2}}
\]
and that (using $G(x)=F(x)-x^{-\alpha }H(x)$)
\begin{eqnarray*}
\lefteqn{\Phi _{N(t)}(u)} \\
&& \hspace{-7mm} =1+\frac{(e^{iu}-1)}{\overline{G}(t)}\; \frac{1}{(1-\frac{
(e^{iu}-1)}{\overline{G}(t)}G(t))^{2}}\; \frac{t^{-\alpha }H(t)+G(t)
\overline{G}(t)-(e^{iu}-1)G^{2}(t)}{\overline{G}(t)} \\
&& \hspace{-7mm} =1+\frac{(e^{iu}-1)u}{u\overline{G}(t)}\; \frac{1}{(1-
\frac{(e^{iu}-1)u}{u\overline{G}(t)}G(t))^{2}} \left[\frac{t^{-\alpha }H(t)
}{\overline{G}(t)}+G(t)-\frac{(e^{iu}-1)uG^{2}(t)}{u\overline{G}(t)}\right].
\end{eqnarray*}
Now we use the approximation $(e^{iu}-1)/u\approx i$ for $u \approx 0$ and replace $u$ by $
v\overline{G}(t)$ (which converges to $0$ as $t\rightarrow \infty $). We find that for $t \rightarrow \infty$ 
\begin{eqnarray*}
\Phi _{N(t)}(v\overline{G}(t)) &\rightarrow &1+iv \; \frac{1}{(1-iv)^{2}}
\; \left(\frac{\alpha -\theta }{\alpha }+1-iv\right) \\
&=& \frac{1-\theta /\alpha }{(1-iv)^{2}}+\frac{\theta /\alpha }{1-iv}.
\end{eqnarray*}
The second result (ii) follows from (i).
\qed

\vspace{2mm}
Notice that for the classical Poisson process we have that $\frac{N(t)}{t}$ converges to a constant random variable $\lambda$.
\vspace{2mm}

\begin{rem}
\rm {If $m_{\alpha} = E(T_1^{\alpha})< \infty$, then Theorem 2 shows that we have}
$$
\lim\limits_{t\to \infty} \frac{R(t)}{t^{\alpha}}= \frac{2}{m_\alpha}, \quad
\lim_{t \to \infty} \frac{\mathbf{E}N^2(t)}{t^{2\alpha}} = \frac{6}{m_\alpha^2}, \quad
\lim_{t \to \infty} \frac{Var(N(t))}{t^{2\alpha}} = \frac{2}{m_\alpha^2}.
$$
\end{rem}

In the next result we discuss the sum $S_{n}$. For this summation process,
we have the following identity:
\[
P(S_{n}\leq t)=P(N(t)>n)\text{.}
\]
For convenience we set $Q(t)=1/\overline{G}(t)$.
\begin{prop}
Suppose that $H\in RV_{\theta }$ where $0\leq \theta <\alpha $. Then there
exists an increasing function $U(x)$ so that $U(Q(x))\sim x$ and $
S_n/U(n)\Longrightarrow Z^{-1/(\alpha -\theta )}$, where $Z$ is given in
the previous result.
\end{prop}

\noindent{\bf Proof.}
Using $Q(x)=1/\overline{G}(x)\in RV_{\alpha -\theta }$ we have $
N(t)/Q(t)\Longrightarrow Z$. Since $\alpha -\theta >0$, $Q(x)$ is
asymptotically equal to a strictly increasing function $V(x)\in RV_{\alpha
-\theta }$ (see \cite {BGT}, Section 1.5.2, Theorem 1.5.4, p.23) and $
N(t)/V(t)\Longrightarrow Z$ in the sense of distribution. We denote the inverse of $V(x)$ by $U(x)$ and
then $U\in RV_{1/(\alpha -\theta )}$. Clearly $N(t)/V(t)\Longrightarrow Z$
implies that $N(U(t))/t\Longrightarrow Z$. Now we have
$$
\mathbf{P} \left\{S_{n}\leq U(n)t\right\}=\mathbf{P} \left\{N(U(n)t)>n\right\}= \mathbf{P}  \left\{\frac{N(U(n)t)}{V(U(n)t)}> 
\frac{n}{ V(U(n)t)}\right\}.
$$
Since $V(U(n)t)/n\sim V(U(n))t^{\alpha -\theta }/n\rightarrow t^{\alpha
-\theta }$, we have 
$$
\mathbf{P} \left\{S_n \leq U(n)t \right\} \rightarrow \mathbf{P} \left\{ Z\geq t^{\theta - \alpha } \right\} =\mathbf{P} \left\{Z^{-1/(\alpha -\theta )}\leq t\right\}.
$$
We conclude that $S_n/U(n)\Longrightarrow Z^{-1/(\alpha -\theta )}$ in the sense of distribution. The
construction shows that $U(Q(x))\sim U(V(x))\sim x$.

\vspace{2mm}

\begin{Ex}\label{ex:RV_Pareto}
Let $f(x)=\beta x^{-\beta -1},x\geq 1$ denote the Pareto density with
parameter $\beta >0$. Clearly we have $\overline{F}(x)=x^{-\beta },x\geq 1$.
For $H(x)$ we find that
\[
H(x)=\int_{1}^{x}y^{\alpha }f(y)dy=\beta \int_{1}^{x}y^{\alpha -\beta -1}dy
\text{.}
\]
We have the following cases.

1) If $\alpha <\beta $, then $H(x)\rightarrow m(\alpha )<\infty $;

2) If $\alpha =\beta $, then $H(x)=\alpha \log x$ and $H\in RV_{0}$ with $
xH^{\prime }(x)/H(x)\rightarrow 0$.

In this case we have $\overline{G}(x)\sim \alpha x^{-\alpha }\log x$ so that
we may take $V(x)=x^{\alpha }/(\alpha \log x)$. To find its inverse
function, we set
\[
\frac{x^{\alpha }}{\alpha \log x}=y\text{.}
\]
We find that $x=(\alpha y\log x)^{1/\alpha }$ and that $\alpha \log x=(\log
\alpha +\log y+\log \log x)$. It follows that $\alpha \log x\sim \log y$ and
\[
U(y)\sim (y\log y)^{1/\alpha }\text{.}
\]

For this example we find

\begin{itemize}
\item $\alpha x^{-\alpha }\log xR(x)\rightarrow 2$;

\item $x^{1-\alpha }\log x(R(x+h)-R(x))\rightarrow 2h$;

\item $\alpha t^{-\alpha }\log tN(t)\Longrightarrow Z$ and $\Phi
_{Z}(v)=(1-iv)^{-2}$;

\item $S_{n}/(n\log n)^{1/\alpha }\Longrightarrow Z^{-1/\alpha }$.
\end{itemize}

3) If $\alpha >\beta $, then $H(x)=\beta (x^{\alpha -\beta }-1)/(\alpha
-\beta )\sim \beta x^{\alpha -\beta }/(\alpha -\beta )$ and $H\in RV_{\theta
}$ with $\theta =\alpha -\beta >0$

In this case we have $\overline{G}(x)\sim \alpha x^{-\beta }/\theta $. We
may take $V(x)=\theta x^{\beta }/\alpha $ and find its inverse $U(y)=(\alpha
y/\theta )^{1/\beta }$.

For this example we find

\begin{itemize}
\item $x^{-\beta }R(x)\rightarrow (\alpha -\beta )(\alpha +\beta )/\alpha
^{2}$ \

\item $x^{1-\beta }(R(x+h)-R(x))\rightarrow eh$ with $e=\theta \beta (\alpha
+\theta )/\alpha ^{2}$

\item $\alpha t^{-\beta }N(t)/\theta \Longrightarrow Z$

\item $S_{n}/(\alpha n/\theta )^{1/\beta }\Longrightarrow Z^{-1/\beta }$
\end{itemize}
\end{Ex}
\bigskip

\begin{Ex}\label{ex:RV_Cauchy} 
For the standard Cauchy density we have $f(x)=\frac{2}{\pi (1+x^{2})}$,  $x\geq 0$.
Clearly we have $f(x)\sim 2 \pi^{-1} x^{-2}$. We can use same analysis as in the
previous example with $\beta =1$.
\end{Ex}
\begin{Ex}\label{ex:RV_student} 
The one-sided $t$-distribution has a density of the form $f(x)\sim c(\beta
)x^{-\beta -1}$ as $x\rightarrow \infty $. We can use the same analysis as
in example \ref{ex:RV_Pareto}.
\end{Ex}
{\bf Acknowledgements.} This paper is a part of project "First order Kendall maximal autoregressive processes and their applications", which is carried out within the POWROTY/REINTEGRATION programme of the Foundation for Polish Science co-financed by the European Union under the European Regional Development Fund.
%%%%%%%%%%%%%%%%%%%%%%%%%%%%%%%%%%%%%%%

\end{document}